\documentclass[12pt,a4paper]{article}

\usepackage{graphicx}				
\usepackage{color}				
\usepackage{amsmath}				
\usepackage{amsthm}
\usepackage{amssymb}
\usepackage{theoremref}
\usepackage{amsfonts}
\usepackage{here}
\usepackage[utf8x]{inputenc}
\usepackage{pdfpages}

\newtheorem{theorem}{Theorem}
\newtheorem{lemma}{Lemma}

\newtheorem{corollary}{Corollary}

\newcommand{\seg}{\textrm{segm}_+}

\DeclareMathOperator{\dist}{dist}

\author{Paul F.X. Müller and Katharina Riegler}
\date{ \today}
\title{Radial Variation of Positive Harmonic Functions on Lipschitz Domains}

\begin{document}
\maketitle
\begin{abstract}
In this paper we present the necessary modifications of the proof by Mozolyako and Havin in \cite{MozHav} to work for Lipschitz domains instead of $\mathcal{C}^2$ domains. 
\end{abstract}

\section{Introduction}

In this article we will present a result (\thref{Lipschitz}) on the variation of positive harmonic functions on Lipschitz domains. The development of results of this type starts with a theorem by Bourgain (\cite{Bourgain}) asserting that the radial variation of a positive harmonic function on the unit disc in $\mathbb{C}$ is bounded at least in one direction.

In 2016, Mozolyako and Havin published the following result for subdomains of $\mathbb{R}^{n+1}$ with a $\mathcal{C}^2$ boundary.
\begin{theorem}[Mozolyako $\&$ Havin, 2016]
\thlabel{MozHav}
Let $u$ be a positive harmonic function on a domain $O\subset\mathbb{R}^{n+1}$ with a $\mathcal{C}^2$ boundary with and an fixed inner point $z_0$. Let $N(x)$ denoted the vector which is normal to the boundary at the point $x\in \partial O$ pointing inside the domain. Let $r$ be a positive function on the boundary $\partial O$ such that $(x,x+r(x)N(x)]\subset O$ for all $x$. Then for all surface balls $E$ with $\omega^{z_0}(E,O)\geq c$ there is a $x\in E$ such that
\begin{equation*}
\int\limits_{0}\limits^{r(x)} |\nabla u(x+yN(x))| dy<cu(z_0),
\end{equation*}
where the constant $c_1$ only depends on the $\mathcal{C}^2$-constant of the boundary of the domain, the constant $c$ and the Harnack distance between $z_0$ and $x+r(x)N(x)$. 
\end{theorem}

On the first page of \cite{MozHav}, Havin and Mozolyako state that the $\mathcal{C}^2$ condition can be relaxed considerably. Motivated by our work on radial variation on Bloch functions in $\mathbb{R}^d$ we were especially interested in the Havin Mozolyako result for Lipschitz domains. Here it is important to note that Havin and Mozolyako point out that their result extends well beyond the $\mathcal{C}^2$ class of domains (mentioning harmonic functions on Riemannian manifolds etc). In the present paper we prove \thref{MozHav} holds true for Lipschitz domains. Thus, the Havin Mozolyako Theorem in the version of \thref{Lipschitz} provides the starting point for our work in \cite{eigenes}.

\begin{theorem}
\thlabel{Lipschitz}
Let $u$ be a positive harmonic function on a Lipschitz domain $O\subset\mathbb{R}^d$ with starcenter $z_0$ and boundary $D$. Let $N(p)$ be a direction "well-inside the domain at $p$" and $r$ a positive function on $D$ such that $[p,p+r(p)N(p)]\subset O$ for all $p\in D$. Then for all surface balls $E\subset D$ with $\omega^{z_0}(E,O)\geq c$ there is a $p_0\in E$ and a harmonic majorant $H$ of the gradient such that 
\begin{equation*}
\int\limits_{0}\limits^{r(p_0)}  H(p+yN(p)) dy < c_1u(z_0)
\end{equation*}
where the constant $c_1$ only depends on the Lipschitz constant of the domain, the constant $c$ and the Harnack distance between $z_0$ and $p_0+r(p_0)N(p_0)$. 
\end{theorem}
A surface ball is a the intersection of a ball centered on the boundary and the boundary. 
We call a direction "well-inside the domain at $p$" if it is within a cone with apex $p$ inside the domain such that also the cone with an opening angle twice as large is contained in the domain.

The article is dedicated to give a proof of \thref{Lipschitz}. It is along the lines of the proof of Mozolyako and Havin in \cite{MozHav}. Therefore we want to point out the similarities and differences before going into the details of the proof. 
First we isolate the places in \cite{MozHav} where the $\mathcal{C}^2$ assumption is used explicitly:
\begin{itemize}
\item 
On page 5 (Proof of (3.7)) the authors select a ball tangent to the boundary. They use the $\mathcal{C}^2$ assumption to estimate the Poisson kernels of a domain using Poisson kernels of balls (\cite{MozHav} Corollary 1, page 27). Havin and Mozolyako used Poisson kernels with respect to surface measure. The approximation of the domain allowed Havin and Mozolyako to prove (3.7) in \cite{MozHav} which is 
\begin{equation}
\label{(3.7)}
\frac{p_{y_2}}{p_{y_1}}\leq c(S)\frac{y_2}{y_1} 
\end{equation} for $y_1\leq y_2$. The constant $c(S)$ depends expressly on the $\mathcal{C}^2$-norm of the parametrization of the domain. The inequality is systematically exploited throughout \cite{MozHav}. Emphasis is on the exponent $1$ on the right-hand side of (\ref{(3.7)}).
\item Continuity of the Poisson kernel is used throughout the paper. 
\end{itemize}

For Lipschitz domains the above approximation of the domain is not possible. The resulting kernel estimates do not hold and continuity is not guaranteed. Therefore some modifications are necessary: 

\begin{itemize}
\item We replace the Poisson kernel with respect to surface measure by the Martin kernel with respect to harmonic measure with pole at a fixed point. Correspondingly we replace surface measures in the definition of the kernel operators $K_y,C_y,B_y, \Omega_y, \cdots$ by harmonic measures with pole at a fixed point. 

\item 
For Martin kernels we obtain the following substitute for (\ref{(3.7)}):
\begin{equation}
\label{austausch}
\frac{k_{y_2}}{k_{y_1}}\leq c(S)\left(\frac{y_2}{y_1}\right)^\alpha
\end{equation}
for $y_1\leq y_2$, where $\alpha$ is a large number.
\end{itemize}

In the following sections we execute the basic line of proof of \cite{MozHav} using (\ref{austausch}) instead of (\ref{(3.7)}). This led to some regrouping of the original Havin and Mozolyako argument, but no serious obstacle arose.  

\begin{itemize}
\item We establish convergence of the $\Pi^\mu$ kernels giving rise to the central kernel $\omega_\Delta$  at the end of section \ref{construction}.
\item We give a proof of \thref{omega_delta}, which states properties corresponding to (1)-(5) on page 22 in \cite{MozHav}.
\item Finally, we note that Havin and Mozolyako obtained the differential equation, $\Phi$-Property and the other Lemmata at the end of the proof using (1)-(5) on page 22 in \cite{MozHav} only. Therefore those proofs carry over easily, as is seen towards the end of the proof given here.   
\end{itemize}

\section{Preliminaries}
\label{prel}

\subsection{Notation}
We use the notation $B$ or $B^d$ for the unit ball of $\mathbb{R}^d$, $S$ for its boundary, the unit sphere, $\mathbb{D}$ for the unit disc of the complex plane $\mathbb{C}$, $B_r$ or $B(r)$ for the ball with center $0$ and radius $r$. The Euclidean distance between two points or a set and a point will be denoted by $d(\cdot,\cdot)$ and the diameter of a set $A$ with $\textrm{diam}(A)$.
For domains $E$ their Euclidean boundary is denoted by $\partial E$, the inward unit vector of a point $x$ of $\partial E$, if it is well-defined, by $N(x)$.

\subsection{Poisson Kernel}
The Poisson kernel on the unit ball $p:B^d\times S \rightarrow\mathbb{R}$ is given by $p(z,\zeta):=\frac{1-|z|^2}{\omega_{d-1}|\zeta-z|^d}$ for $|z|<1$, $|\zeta|=1$ and $\omega_{d-1}$ the surface area of the unit sphere.

\subsection{Harnack's Inequality}
We will use Harnack's inequality to compare values of positive harmonic functions and to get a bound for their gradients. 
\begin{theorem}
Let $h$ be a positive harmonic function on a ball ${B}(x_0,R)\subset\mathbb{R}^d$ then for all $x$ such that 
$\dist(x,x_0)=r<R$ we have
\[\frac{(1-\frac{r}{R})}{(1+\frac{r}{R})^{d-1}} f(x_0)\leq f(x)\leq \frac{(1+\frac{r}{R})}{(1-\frac{r}{R})^{d-1}} f(x_0)\]
\end{theorem}

\begin{corollary}
Let $f$ be a positive harmonic function on a domain $\Omega$ and $x\in \Omega$ then we have the following estimate for the gradient:
\[|\nabla f|(x)\leq \frac{f(x)}{\dist(x,\partial \Omega)}\]
\end{corollary}

See \cite[Theorem 1.3.1]{Ransford} and \cite[Section 2.4]{Helms}.

\subsection{Harmonic Measure}
Given a set $E$, we use the notation $w^{z_0}(F,E)$ for the harmonic measure with pole $z_0\in E$ of $F\subset \partial E$. It is a probability measure on $\partial E$. In $z$ it is a harmonic function on $E$ solving the Dirichlet problem with boundary data equal to the indicator function of $F$. And it is the probability of the Brownian motion started in $z_0$, stopped at $\partial E$, to be stopped in $F$.
See \cite{Ransford}. 

\subsection{Green's Function}
A Green's function for a domain $\Omega\subset\mathbb{C}$ is a function $g:\Omega\times\Omega\rightarrow (-\infty,\infty]$ such that for each $w\in\Omega$
\begin{enumerate}
\item $g(\cdot,w)$ is harmonic on $\Omega\setminus \{w\}$ and bounded outside each neighbourhood of $w$
\item $g(w,w)=\infty$ and as $z\rightarrow w$
\begin{equation*}
g(z,w)=
\begin{cases}
\log|z|+O(1)&w=\infty\\
-\log|z-w|+O(1)&w\neq \infty
\end{cases}
\end{equation*}
\item $g(z,w)\rightarrow 0$ as $z\rightarrow\zeta$ and $\zeta \in \partial\Omega$.
\end{enumerate}
In $\mathbb{R}^d$ the singularity is of the form $\frac{1}{|z-w|^{d-2}}$. See \cite[Section 4.4]{Ransford}.

\subsection{Lipschitz domains}
We will prove and use a theorem on positive harmonic functions on a Lipschitz domain. Therefore we want to give the definition of a Lipschitz domain. 
A domain $\Omega\subset \mathbb{R}^d$ is called Lipschitz domain if for every point $p\in\partial \Omega$ there is, up to translation an rotation, a Lipschitz function $g:\mathbb{R}^{d-1}\rightarrow\mathbb{R}$ and real numers $r>0$, $h>0$ such that
\begin{itemize}
\item $\Omega\cap \mathcal{C}=\{(x,y)\in\mathbb{R}^{d-1}\times\mathbb{R}: |x|<r, g(x)<y<h\}$
\item $\partial \Omega \cap \mathcal{C}=\{(x,y)\in\mathbb{R}^{d-1}\times\mathbb{R}: |x|<r, g(x)=y\} $,
\end{itemize}
where
\begin{equation}
\mathcal{C}:=\{(x,y)\in\mathbb{R}^{d-1}\times\mathbb{R}: |x|<r, -h<y<h\} .
\end{equation}

\subsection{Martin Boundary and Martin Kernel}
We will make extensive use of the Martin kernel of a (Lipschitz) domain $\Omega$ denoted by $k^{\Omega}$. It will, for example, substitute the Poisson kernel used in \cite{MozHav}.

To define the Martin boundary of a domain we consider $M(x,y):= \frac{g(x,y)}{g(x_0,y)}$ where $x_0$ is a fixed point in the domain. The function $x\mapsto M(x,y)$ is continuous for $y\in\Omega\setminus\{x\}$. We now use the theorem of Constantinescu-Cornea (see \cite[Theorem 7.2]{Bass} or \cite[Theorem 12.1]{Helms}) to get a compact set $\Omega^*$, unique up to homeomorphisms such that
\begin{enumerate}
\item $\Omega$ is a dense subset of $\Omega^*$,
\item for each $y\in\Omega$ the function $x\mapsto M(x,y)$ has a continuous extension to $\Omega^*$ and
\item the extended functions separate points of $\Omega^*\setminus\Omega$. 
\end{enumerate}
The set $\Omega^*\setminus\Omega$ is the Martin boundary of $\Omega$ and denoted by $\partial_M \Omega$. The extensions of $M$ are called Martin kernels and denoted by $k^{\Omega}$. Martin kernels provide the following fundamental representation theorem for positive harmonic functions. 
\begin{theorem}
\thlabel{rep theorem}
For every positive harmonic function $h$ on $\Omega$ there is a measure $\nu$ concentrated on $\partial_M\Omega$ such that 
\begin{equation*}
h(x)=\int k^\Omega(x,y)d\nu(y).
\end{equation*}
\end{theorem}
\paragraph*{Remarks}
\begin{itemize}
\item In the case of Lipschitz domains the Martin boundary and the Euclidean boundary coincide.
\item For any $\zeta\in \partial\Omega$ the function $x\mapsto k(x,\zeta)$ is a positive harmonic function.
\item The Martin kernel is continuous on $\Omega\times \partial \Omega$. 
\end{itemize}
See \cite[Section II.7]{Bass} or \cite[Chapter 12]{Helms} and \cite{HuntW}.

\section{Proof of Theorem 2}
\subsection{Definition of domain}
Let $\Phi:\mathbb{R}^{d-1}\rightarrow\mathbb{R}$ a Lipschitz function satisfying the following conditions:
\begin{itemize}
\item $\Phi(0)=0$,
\item there exists $r\in (0,1)$ such that for any $x\in\mathbb{R}^{d-1}\setminus B^{d-1}(0,r)$ we have $\Phi(x)=0$.
\end{itemize}
The the near half space $O$ is given by $O:=\{(x,y)\in \mathbb{R}^d:x\in\mathbb{R}^{d-1}, y>\Phi(x)\}$. Its boundary $S$ is the graph of $\Phi$.

From now on $u$ is a fixed positive harmonic function with \[\lim\limits_{(x^*,y)\rightarrow(x_0^*,0)} u(x^*,y)=\lim\limits_{|x^*|+y\rightarrow\infty} u(x^*,y)=0\] for $x_0^*\in \mathbb{R}^{d-1}$ with $|x_0^*|>r$.

Throughout the chapter we use the notation $x_y=x+y \vec{e}_d$, where $x\in \mathbb{R}^d$ $y\in\mathbb{R}$ and $\vec{e}_d$ is the last of the standard basis vectors in $\mathbb{R}^d$. Analogously, the sets are shifted: $E_y:= E+y\vec{e}_d$.  For functions $\phi$ on subsets of $\mathbb{R}^d$ we write $\phi_y(x)=\phi(x_y)$. We also fix a point $z_0$ that will be the pole of the harmonic measure that we use.  

As $\Phi$ is a Lipschitz function we know that 
\begin{equation}
\label{distance comparable}
y\geq \dist(x_y,S)\geq c(S)y
\end{equation}
for all $x\in\S$ and $y>0$.

\subsection{Definition and Properties of basic kernels}
In the following by a kernel we mean a function defined on $S\times S$. The kernels are denoted by lowercase letters and the corresponding integral operators by the respective uppercase letter, for example $Qu:=\int\limits_{S}q(x,\zeta)u(\zeta)d\omega^{z_0}(\zeta)$, where $\omega^{z_0}$ is the harmonic measure on $O$ with pole at $z_0$.
The composition of two kernels is defined as
$(p\circ q )(x,\zeta):=\int \limits_{S}p(x,\eta)q(\eta,\zeta)d\omega^{z_0}(\eta)$.

We now give the definition of the three basic kernels in the proof. The first familiy of kernels is given by \[k_y(x,\zeta):=k^O(x_y,\zeta),\] where $k$ is the Martin kernel of $O$. 
The other two families also depend on the positive harmonic function $u$. We have \[c_y(x,\zeta):=\frac{\partial^1 k}{\partial \sigma(x_{2y})}(x_y,\zeta),\] where $\sigma(x):=\textrm{sgn}((\nabla u)(x))$ and $\textrm{sgn}(q)$ is the normalization of $q$ ($\textrm{sgn}(q)=0$ for $q=0$, and then also $\frac{\partial^1 k}{\partial \sigma(q)}=0$). The last kernel is given by 
\begin{equation}
\label{def b_delta}
b_y:=k_y\circ c_y.
\end{equation}

\paragraph*{Properties of $k_y$:}
\begin{lemma}
\thlabel{prop ky}
The kernel $k_y$ has the following properties:
\begin{enumerate}
\item \label{ky 0} For any positive harmonic function $u$ on $O$ we have for all $x\in S$ that $K_{y_2}(u_{y_1}\big|_S)(x)=u_{y_1+y_2}(x)$
\item \label{ky 1} For any $y_1,y_2>0$ we have the following semi-group property: $k_{y_1+y_2}=k_{y_1}\circ k_{y_2}$
\item \label{ky 2} For $y_1\leq y_2$ we have $\frac{k_{y_2}}{k_{y_1}}\leq c\left(\frac{y_2}{y_1}\right)^{\alpha}$, where $\alpha=c(S)$ is a constant only depending on $S$.
\item \label{ky 3} $K_y(1)=1$
\end{enumerate}
\begin{proof}~\\
Proof of \ref{ky 0}: \\By definition $K_{y_2}(u_{y_1}\big|_S)(x)$ is equal to $\int\limits_S k_{y_2}(x,\zeta)u_{y_1}\big|_S(\zeta) d\omega^{z_0}(\zeta)(x)$, which is the harmonic continuation of $u_{y_1}\big|_S$ evalutated at $x_{y_2}$. This is $u_{y_1}(x_{y_2})=u(x_{y_1+y_2})$.\\
Proof of \ref{ky 1}: \\This is a special case of \ref{ky 0}.\\
Proof of \ref{ky 2}:\\
Let $c=c(S)$ be a constant such that $\dist(x_y,S)\geq cy$ as in (\ref{distance comparable}). We then know that for any $x\in S$ and any $y>0$ the ball $B(x_y, cy)$ is contained in $O$. Therefore $k_y(\cdot, \zeta)$ is a positive harmonic function on $B(x_y, cy)$. We will thus be able to apply Harnack's inequality in this setting. 
Let $n\in\mathbb{N}$ be such that $y_2 (1-\frac{c}{2})^n\leq y_1 \leq y_2 (1-\frac{c}{2})^{n-1}$.
Let $a_1=y_2$ and $a_k=a_1(1-\frac{c}{2})^k$ for $k=1,...,n$. We then know by Harnack's inequality that
\begin{equation*}
k(x_{y_2},\zeta)=k(a_1,\zeta)\leq c^{n-1}k(x_{a_n},\zeta)\leq c^nk(x_{y_1},\zeta),
\end{equation*}
where $c=c(d)$ only depends on the dimension of the domain. 
As $(1-\frac{c}{2})^n\leq \frac{y_1}{y_2}\leq (1-\frac{c}{2})^{n-1}$ we have that $c^n$ is of the form $\frac{y_2}{y_1}^{c(S)}$ for a constant $c(S)$. \\
Proof of \ref{ky 3}:\\By definition we have for all $x\in S$ and all $y>0$:
\[K_y(1)=\int\limits_{S}k(x_y,\zeta)d\zeta=1.\]
\end{proof}
\end{lemma}


\paragraph*{Properties of $c_y$:}
\begin{lemma}
The kernel $c_y$ has the following properties:
\begin{enumerate}
\item \label{cy 1}
$|\nabla u(x_{2y})|=C_y(u_y)(x)$
\item \label{cy 2}
$|c_y(x,\zeta)|\leq c(S)\frac{k_y(x,\zeta)}{y}$
\item \label{cy 3}
$C_y(1)=0$ 
\end{enumerate}
\begin{proof}~\\
Proof of \ref{cy 1}:\\
We use the definition of $c_y$ to obtain:
\begin{align*}
|\nabla u(x_{2y})|&=\langle \nabla u (x_{2y}),\sigma (x_{2y})\rangle=\\
&=\langle\int\limits_{S}\nabla^1 k_y(x ,\zeta)u_y(\zeta)d\omega^{z_0}(\zeta),\sigma(x_{2y})\rangle=C_y(u_y).
\end{align*}
Proof of \ref{cy 2}:\\ By definition of the kernel, the Harnack inequality and the fact that
$\dist (x_y,S) \geq c(S)y$ we get:
\begin{align*}
c_y(x,\zeta)|&=|\frac{\partial^1 k}{\partial \sigma(x_{2y})}(x_y,\zeta)|\leq\\
&\leq |\nabla ^1k(x_y,\zeta)|\leq c\frac{k(x_y,\zeta)}{\dist(x_y,S)} \leq c(S)\frac{k(x_y,\zeta)}{y}
\end{align*}
Proof of \ref{cy 3}:\\
We differentiate property \ref{ky 3} of $k_y$ to obtain $C_y(1)=0.$

\end{proof}
\end{lemma}


\paragraph*{Properties of $b_y$:}
\begin{lemma}
\thlabel{prop by}
The kernel $b_y$ has the following properties:
\begin{enumerate}
\item \label{by 1} $|b_y|\leq c(S)\frac{k_y}{y}$
\item \label{by 2} $B_y(1)=0$
\item \label{by 3} $(x,\zeta)\mapsto b_y(x,\zeta)$ is continuous on $S\times S$
\end{enumerate}
\begin{proof}
Proof of \ref{by 1}:\\
By definition of $b_y$, the positivity of $k_y$, property \ref{cy 2} of $c_y$ and properties \ref{ky 1} and \ref{ky 2} of $k_y$ we have:
\[|b_y|=|k_y\circ c_y|\leq c(S)\frac{k_y\circ k_y}{y}\leq c(S)\frac{k_y}{y}.\]
Proof of \ref{by 2}: \\
Using property \ref{cy 3} of $c_y$ and get $B_y(1)=K_y(C_y(1))=K_y(0)=0$.\\
Proof of \ref{by 3}: \\
By definition $b_y(x,\zeta)=(k_y\circ c_y) (x,\zeta)=\int\limits_{S}k_y(x,\eta)c_y(\eta,\zeta)d\omega^{z_0}(\eta)$.
For a fixed $y$ the kernel $k_y$ is bounded and therefore the constant function is an integrable majorant of $k_y(x,\eta)c_y(\eta,\zeta)$ (we use property (\ref{cy 2}) of $c_y$). 

As $c_y(x,\zeta)=\frac{\partial^1 k}{\partial \sigma(x_{2y})}(x_y,\zeta)$, it is continuous except for the zeros of $\nabla u$. The harmonicity of $u$ yields analyticity of $\nabla u$ and therefore we only have a discrete set of zeros for any $y>0$. Therefore $b_y$ is continuous. 

\end{proof}
\end{lemma}

\subsection{Variations}
We will use the following mean vertical variation
\[V(x):=\int\limits_0\limits^1 B_y(u_y)dy\]
for any point $x\in S$. 

By definition of the operator $B_y$ we have $V(x)=\int\limits_{0}\limits^{1} K_y(C_y(u_y)) dy=\int\limits_{0}\limits^{1} K_y(|\nabla u(\cdot_{2y})|) dy\geq \int\limits_{0}\limits^{1} |\nabla u(x_{3y})|dy$. Therefore the mean vertical variation is (up to the constant $\frac{1}{3}$) greater than the variation along the vertical line. 

\subsection{The Main Lemma}
In this section we state the central lemma. The proof will be done in Section \ref{Proof of Lemma}. 
\begin{lemma}
\thlabel{main lemma}
For any ball $\mathbb{B}$ centered on the boundary $S$ of the almost-half-space $O$ there is a point $x$ in $\mathbb{B}\cap S$ such that the variation $\int\limits_{0}\limits^{1} B_y(u_y)(x) dy$ of the positive harmonic function $u$ is bounded by $cu(z_1)$, where $z_1$ is a fixed point with $z_1=x_y$ where $x\in S$ and $y>1$. The constant $c$ depends only on the Lipschitz constant of the function defining $O$ and the radius of the ball $\mathbb{B}$. 
\end{lemma}

\subsection{Construction of kernel $\omega_\Delta$}
\label{construction}
\paragraph{Notation for segments:}
The set of all non-degenerate compact intervals (segments) in $(0,\infty)$ will be denoted by $\seg$. For $\Delta\in\seg$ we use the following notation for the minimum, maximum and length of $\Delta$:
\begin{align*}
m(\Delta)&:=\min(\Delta)\\
M(\Delta)&:=\max(\Delta)\\
|\Delta|&:=M(\Delta)-m(\Delta).
\end{align*}

\paragraph{Notation for Partitions of a Segment}
We want to partition segments $\Delta$ into smaller segments. Therefore we call a finite set $\mu\subset\seg$ a partition of $\Delta$ if it is a set of non-overlapping intervals, the union of which is $\Delta$. The elements of the partition will be denoted by $j_s$ for $s=1,...,K$. We will number them such that $0< m(\Delta )=m(j_1)<M(j_1)=m(j_2)<...<m(j_K)<M(j_K)=M(\Delta )$. The mesh of a partition $\mu$ is the length of the longest intervel i.e. $\lambda(\mu):=\max\limits_{j\in\mu}|j|$. A partition is called regular if $\lambda(\mu)\leq \frac{2|\Delta|}{K}$.
For refinements of partitions we write $\mu_2\succ \mu_1$ if any element of $\mu_2$ lies in some element of $\mu_1$.

The dyadic partition of $\Delta$ consisting of the segments $\Delta\cap [\frac{s}{2^n},\frac{s+1}{2^n}]$, $s\in\mathbb{Z}$ will be denoted by $\tau_n(\Delta)$. The dyadic partition is regular and $\tau_{n+1}(\Delta)\succ\tau_n(\Delta)$.
\paragraph{Definition and properties of $b_\Delta$:}
We introduce a new kernel $b_\Delta$:

\[b_\Delta(x,\zeta):=\int \limits_{\Delta} b_y(x,\zeta) dy\]
\begin{lemma}
The kernel $b_\Delta$ has the following properties:
\begin{enumerate}
\item \label{bd 1}$|b_\Delta|\leq c(S)\frac{M(\Delta)^{\alpha-1}}{m(\Delta)^\alpha}|\Delta|k_{m(\Delta)}$
\item \label{bd 2} $b_\Delta$ is continuous on $S\times S$
\end{enumerate}
\begin{proof}
Proof of \ref{bd 1}:\\
Using property (1) of $b_y$ we get 
\begin{equation*}
|b_\Delta|\leq \int\limits_{\Delta}c(S)\frac{k_y}{y}dy\leq \int\limits_{\Delta} c(S) k_{m(\Delta)} \frac{y^{\alpha-1}}{m(\Delta)^\alpha} dy\leq c(S) |\Delta|\frac{M(\Delta)^{\alpha-1}}{m(\Delta)^\alpha} k_{m(\Delta)}
\end{equation*}
Proof of \ref{bd 2}: We use the continuity of $b_y$ and the constant function as an integrable majorant of $|b_y(x,\zeta)|\leq c(S)\frac{k_y(x,\zeta)}{y}$ to obtain \[\lim\limits_{(x,\zeta)\rightarrow(x_0,\zeta_0)}b_{\Delta}(x,\zeta)=\lim\limits_{(x,\zeta)\rightarrow(x_0,\zeta_0)}\int\limits_{\Delta}b_y(x,\zeta)dy=b_\Delta(x_0,\zeta_0).\]
\end{proof}
\end{lemma}

\paragraph{Definition and properties of $\tilde{\omega}_\Delta$}
We introduce a new kernel $\tilde{\omega}_\Delta$ which is dependent on a small real number $\epsilon$. The value of $\epsilon$ will be fixed at a later point
\[\tilde{\omega}_\Delta:=k_{|\Delta|}-\epsilon b_\Delta.\]\\

\paragraph{Definition and Properties of $\Pi^\mu$}
Given a partition $\mu$ of the segment $\Delta$ we define the kernel $\Pi^\mu$ by
\[\Pi^\mu:=\tilde{\omega}_{j_K}\circ\tilde{\omega}_{j_{K-1}}\circ\cdots\circ\tilde{\omega}_{j_1}.\]
We will now break up the kernel $\Pi^\mu$ into three different parts, that we will be able to treat separately.
By definition
\[\Pi^\mu=(k_{|j_K|}-\epsilon b_{j_K})\circ (k_{|j_{K-1}|}-\epsilon b_{j_{K-1}})\circ\cdots\circ (k_{|j_1|}-\epsilon b_{j_1}).\]
Using the notation $N_q$ for subsets of $\{1,2,\cdots, K\}$ of cardinality $q$ as well as $\pi_l:=r^l_K\circ r^l_{K-1}\circ \cdots \circ r^l_1$ with 
\begin{equation*}
r^l_s= \left\{
\begin{array}{ll}
-\epsilon b_{j_s}\quad & s\in l \\
k_{|j_s|} & \, s\notin l \\
\end{array}
\right.
\end{equation*}
we obtain 
\begin{equation*}
\Pi^\mu = k_{|\Delta|}+\sum\limits_{q=1}\limits^{K}\sum\limits_{l\in N_q}\pi_l.
\end{equation*}
Now we want to isolate the sum where $q=1$. In this case only one $r^l_s$ is of the form $-\epsilon b_{j_s}$ and therefore we get
\begin{equation*}
\sum\limits_{l\in N_1}\pi_l= -\epsilon\sum\limits_{s=1}\limits^{K}k_{M(\Delta)-M(j_s)} \circ b_{j_s}\circ k_{m(j_s)-m(\Delta)}
\end{equation*}
where $k_0$ is understood as the identity.

Taking into account that $b_\Delta=\sum\limits_{s=1}\limits^{K}b_{j_s}$ we obtain
\begin{equation}
\label{pi}
\Pi^\mu = k_{|\Delta|}
-\epsilon b_\Delta
+\epsilon \sum\limits_{s=1}\limits^{K} b_{j_s} - k_{M(\Delta)-M(j_s)} \circ b_{j_s}\circ k_{m(j_s)-m(\Delta)}
+ \sum\limits_{q=2}\limits^{K}\sum\limits_{l\in N_q}\pi_l.
\end{equation}
For further calculations we will use the notation $v_s= k_{M(\Delta)-M(j_s)} \circ b_{j_s}\circ k_{m(j_s)-m(\Delta)}$ and $\rho_\mu=\sum\limits_{q=2}\limits^{K}\sum\limits_{l\in N_q}\pi_l$ and (\ref{pi}) is equivalent to
\begin{equation*}
\Pi^\mu =\tilde{\omega}_\Delta+ \epsilon \sum\limits_{s=1}\limits^{K}( b_{j_s}-v_s) + \rho_\mu
\end{equation*}

We will now show a series of lemmata, that will be of use when we prove the existence of the limit of $\Pi^{\tau_n}$.

\begin{lemma}
We have the following estimate:
\begin{align*}
&\sum\limits_{s=1}\limits^{K} |b_{j_s}-v_s|\leq \\
&\leq c(S)|\Delta|^2 \left(\frac{ M(\Delta)^{\alpha-1}(3M(\Delta)-m(\Delta))^\alpha}{m(\Delta)^{2\alpha+1}}+\frac{(M(\Delta)+|\Delta|)^{\alpha-1} }{m(\Delta)^{\alpha+1}}\right) k_{m(\Delta)}
\end{align*}
\begin{proof}
We will use
\[|b_y|+|c_y|\leq c(S)\frac{k(y)}{y}\]
\[|k_{\theta+\lambda}-k_\theta|\leq c(S)k_\theta \left(\frac{(\theta+\lambda)^\alpha}{\theta^\alpha} -1\right)\]
which we know by the following calculation:
\begin{align*}
|k_{\theta+\lambda}-k_\theta|&\leq \int\limits_0\limits^\lambda |\frac{d}{dt} k_{\theta+t}|dt\\
&\leq c(S)\int\limits_0\limits^\lambda \frac{k_{\theta+t}}{\theta+t} dt\\
&\leq c(S)\int\limits_0\limits^\lambda k_\theta\left(\frac{\theta+t}{\theta}\right)^\alpha\frac{1}{\theta+t}dt\\
&\leq c(S)k_\theta \left(\frac{(\theta+\lambda)^\alpha}{\theta^\alpha} -1\right).
\end{align*}

Now we want to treat $|b_{j_s}-v_s|\leq|b_{j_s}-k_{M(\Delta)-M(j_s)}\circ b_{j_s}|+|k_{M(\Delta)-M(j_s)}\circ b_{j_s}-v_s|=I+II$ in two steps. We start by estimating $I$:

\begin{align*}
I&=|b_{j_s}-k_{M(\Delta)-M(j_s)}\circ b_{j_s}|\\
&\leq |\int\limits_{j_s} b_y dy-k_{M(\Delta)-M(j_s)}\circ\int\limits_{j_s} b_ydy|\\
&\leq \int\limits_{j_s} k_y\circ c_y - k_{M(\Delta)-M(j_s)+y}\circ c_y|dy\\
&\leq \int\limits_{j_s} |k_y-k_{M(\Delta)-M(j_s)+y}|\circ |c_y|dy\\
&\leq c(S)\int\limits_{j_s} \left(\frac{(y+M(\Delta)-M(j_s))^\alpha}{y^\alpha}-1\right) k_y\circ \frac{k_y}{y}dy\\
&\leq c(S)\frac{1}{m(\Delta)} \int\limits_{j_s} \left(\frac{(y+M(\Delta)-M(j_s))^\alpha}{y^\alpha}-1\right) k_{2y}dy\\
&\leq c(S)\frac{1}{m(\Delta)} \int\limits_{j_s} \left(\frac{(y+M(\Delta)-M(j_s))^\alpha}{y^\alpha}-1\right)\frac{(2y)^\alpha}{m(\Delta)^\alpha}  k_{m(\Delta)}dy\\
&\leq c(S)\frac{2^\alpha}{m(\Delta)^{\alpha+1}} k_{m(\Delta)} \int\limits_{j_s} (y+M(\Delta)-M(j_s))^\alpha-y^\alpha dy \\
&\leq c(S)\frac{1}{m(\Delta)^{\alpha+1}} k_{m(\Delta)} \int\limits_{j_s} \int\limits_y\limits^{y+M(\Delta)-M(j_s)}\alpha s^{\alpha-1} ds dy \\
&\leq c(S)\frac{(M(\Delta)+|\Delta|)^{\alpha-1} }{m(\Delta)^{\alpha+1}}  |j_s||\Delta| k_{m(\Delta)} 
\end{align*}
The second part is

\begin{align*}
II&=|k_{M(\Delta)-M(j_s)}\circ b_{j_s}-k_{M(\Delta)-M(j_s)} \circ b_{j_s}\circ k_{m(j_s)-m(\Delta)}| \\
&\leq \int\limits_{j_s}|k_{M(\Delta)-M(j_s)}\circ k_y\circ c_y-k_{M(\Delta)-M(j_s)} \circ k_y\circ c_y\circ k_{m(j_s)-m(\Delta)}|		dy\\
&\leq \int\limits_{j_s}|k_{M(\Delta)-M(j_s)+y}\circ c_y\circ (1- k_{m(j_s)-m(\Delta)})|dy\\
&\leq \int\limits_{j_s}|k_{M(\Delta)-M(j_s)+y}\circ \frac{\partial k_{\frac{y}{2}}}{\partial \sigma(2y)}\circ k_{\frac{y}{2}}\circ (1- 			k_{m(j_s)-m(\Delta)})|dy \\
&\leq c(S)\int\limits_{j_s}k_{M(\Delta)-M(j_s)+y}\circ \frac{2}{y}k_{\frac{y}{2}}\circ |(k_{\frac{y}{2}} - k_{\frac{y}{2}+ 				m(j_s)-m(\Delta)})|dy \\
&\leq c(S) \frac{2}{m(\Delta)}\int\limits_{j_s}k_{M(\Delta)-M(j_s)+2y} \left(\frac{(\frac{y}{2}	+m(j_s)-m(\Delta))^\alpha}				{(\frac{y}{2})^\alpha} -1\right)dy\\
&\leq c(S) \frac{2}{m(\Delta)}\int\limits_{j_s} \frac{(M(\Delta)-M(j_s)+2y)^\alpha}{m(\Delta)^\alpha} k_{m(\Delta)} 					\left(\frac{(\frac{y}{2}+m(j_s)-m(\Delta))^\alpha}{(\frac{y}{2})^\alpha} -1\right)dy\\
&\leq c(S) \frac{2}{m(\Delta)}\frac{(3M(\Delta)-M(j_s))^\alpha}{m(\Delta)^\alpha} k_{m(\Delta)}\int\limits_{j_s} \frac{2^\alpha}		{y^\alpha} \left(		(\frac{y}{2}+m(j_s)-m(\Delta))^\alpha -\frac{y^\alpha}{2^\alpha}\right)dy\\
&\leq c(S) \frac{2}{m(\Delta)}\frac{(3M(\Delta)-M(j_s))^\alpha}{m(\Delta)^\alpha} k_{m(\Delta)}\frac{2^\alpha}{m(\Delta)^\alpha}		\int\limits_{j_s} \int\limits_{\frac{y}{2}}\limits^{\frac{y}{2}+m(j_s)-m(\Delta)} \alpha s^{\alpha-1}dsdy\\
&\leq c(S)|j_s||\Delta| \frac{ M(\Delta)^{\alpha-1}(3M(\Delta)-m(\Delta))^\alpha}{m(\Delta)^{2\alpha+1}} k_{m(\Delta)}.	\\
\end{align*}
Now we have proven that 
\begin{equation*}
|b_{j_s}-v_s|\leq c(S)|j_s||\Delta|\left(\frac{ M(\Delta)^{\alpha-1}(3M(\Delta)-m(\Delta))^\alpha}{m(\Delta)^{2\alpha+1}}+\frac{(M(\Delta)+|\Delta|)^{\alpha-1} }{m(\Delta)^{\alpha+1}}\right) k_{m(\Delta)}.
\end{equation*}
Summing up we obtain 
\begin{equation*}
\sum\limits_{s=1}\limits^{K} |b_{j_s}-v_s|\leq c(S)|\Delta|^2 \left(\frac{ M(\Delta)^{\alpha-1}(3M(\Delta)-m(\Delta))^\alpha}{m(\Delta)^{2\alpha+1}}+\frac{(M(\Delta)+|\Delta|)^{\alpha-1} }{m(\Delta)^{\alpha+1}}\right) k_{m(\Delta)}.
\end{equation*}
\end{proof}
\end{lemma}
The next step is the analysis of $\rho_\mu$.
\begin{lemma}
We have the following estimate for $\rho_\mu$ for regular partitions but without restrictions for $\Delta$:
\[|\rho_\mu|\leq \epsilon^2 \left(\frac{|\Delta|}{m(\Delta)}\right)^2 k_{m(\Delta)}\left(\frac{M(\Delta)+|\Delta|}{m(\Delta)}\right)^\alpha\sum\limits_{q=2}\limits^{\infty}\frac{q^\alpha}{q!}\left(\epsilon \frac{|\Delta|}{m(\Delta)}\right)^{q-2}\left(c(S)\frac{M(\Delta)^{\alpha-1}}{m(\Delta)^{\alpha-1}}\right)^q\]
\begin{proof}

For $l\in N_q$ we use the following notation:
\begin{equation*}
a(l):=\sum\limits_{s\in l}m(j_s)+\sum\limits_{s\notin l}|j_s|.
\end{equation*}
Now we obtain for $l\in N_q$
\begin{align*}
|\pi_l|&\leq k_{a(l)} \prod\limits_{s\in l}\left( \epsilon c(S) \frac{M(\Delta)^{\alpha-1}}{m(\Delta)^{\alpha-1}}\frac{|j_s|}{m(\Delta)} \right)\\
&\leq k_{a(l)}\frac{1}{K^q}\left(\epsilon c(S)\frac{|\Delta|}{m(\Delta)}\frac{M(\Delta)^{\alpha-1}}{m(\Delta)^{\alpha-1}}\right)^q\\
&\leq k_{m(\Delta)}\frac{(qM(\Delta)+|\Delta|)^\alpha}{m(\Delta)^\alpha}\frac{1}{K^q}\left(\epsilon c(S)\frac{|\Delta|}{m(\Delta)}\frac{M(\Delta)^{\alpha-1}}{m(\Delta)^{\alpha-1}}\right)^q.
\end{align*}
Summing up we obtain
\begin{align*}
|\rho_\mu |&\leq \sum\limits_{q=2}\limits^{K}\frac{K^q}{q!}k_{m(\Delta)}\frac{(qM(\Delta)+|\Delta|)^\alpha}{m(\Delta)^\alpha}\frac{1}{K^q}\left(\epsilon c(S)\frac{|\Delta|}{m(\Delta)}\frac{M(\Delta)^{\alpha-1}}{m(\Delta)^{\alpha-1}}\right)^q\\
&\leq \epsilon^2 \left(\frac{|\Delta|}{m(\Delta)}\right)^2 k_{m(\Delta)}\left(\frac{M(\Delta)+|\Delta|}{m(\Delta)}\right)^\alpha\sum\limits_{q=2}\limits^{\infty}\frac{q^\alpha}{q!}\left(\epsilon \frac{|\Delta|}{m(\Delta)}\right)^{q-2}\left(c(S)\frac{M(\Delta)^{\alpha-1}}{m(\Delta)^{\alpha-1}}\right)^q.
\end{align*}
\end{proof}
\end{lemma}

\begin{lemma}
\thlabel{pi_omega}
For regular partitions $\mu$ we have the following estimate:
\begin{align*}
|\Pi^\mu-\tilde{\omega}_\Delta|\leq & \epsilon k_{m(\Delta)} |\Delta|^2\\
&(\epsilon \left(\frac{M(\Delta)+|\Delta|}{m(\Delta)}\right)^\alpha\sum\limits_{q=2}\limits^{\infty}\frac{q^\alpha}{q!}\left(\epsilon\frac{|\Delta|}{m(\Delta)}\right)^{q-2}\left(c(S)\frac{M(\Delta)^{\alpha-1}}{m(\Delta)^{\alpha-1}}\right)^q+\\
&+\epsilon c(S) \left(\frac{ M(\Delta)^{\alpha-1}(3M(\Delta)-m(\Delta))^\alpha}{m(\Delta)^{2\alpha+1}}+\frac{(M(\Delta)+|\Delta|)^{\alpha-1} }{m(\Delta)^{\alpha+1}}\right) ).
\end{align*}
\begin{proof}
The proof consists of collecting information from the lemmata above. 
\end{proof}
\end{lemma}

\begin{lemma}
\thlabel{Pi_mu}
For regular partitions we have the following estimate for $|\Pi^\mu|$

\begin{align*}
|\Pi^\mu|\leq  &k_{|\Delta|}+\epsilon c(S) \frac{M(\Delta)^{\alpha-1}}{m(\Delta)^\alpha}|\Delta|k_{m(\Delta)}\\
&+\epsilon k_{m(\Delta)} |\Delta|^2\\
&(\epsilon \left(\frac{M(\Delta)+|\Delta|}{m(\Delta)}\right)^\alpha\sum\limits_{q=2}\limits^{\infty}\frac{q^\alpha}{q!}\left(\epsilon\frac{|\Delta|}{m(\Delta)}\right)^{q-2}\left(c(S)\frac{M(\Delta)^{\alpha-1}}{m(\Delta)^{\alpha-1}}\right)^q+\\
&+\epsilon c(S) \left(\frac{ M(\Delta)^{\alpha-1}(3M(\Delta)-m(\Delta))^\alpha}{m(\Delta)^{2\alpha+1}}+\frac{(M(\Delta)+|\Delta|)^{\alpha-1} }{m(\Delta)^{\alpha+1}}\right) ).
\end{align*}
\begin{proof}
We obtain the estimate by using \thref{pi_omega} and the definition of $\tilde{\omega}_{\Delta}$.
\end{proof}
\end{lemma}

\begin{lemma}
\thlabel{refinement}
We have the following estimate for the refinement of a partition:
$|\Pi^\tau-\Pi^\sigma|\leq \lambda(\tau)c(\Delta,S)k_{m(\Delta)}$ for $\sigma \succ \tau$.
\begin{proof}
Suppose $\tau=\{\Delta_1,\cdots\Delta_K\}$, $m(\Delta_1)<m(\Delta_2)<\cdots <m(\Delta_K)$. Put $\sigma_k:=\{j\in\sigma:j\subset\Delta_k\}$. Then $\sigma_k$ is a partition of $\Delta_k$ and $\sigma=\bigcup\limits_{k=1}\limits^{K}\sigma_k$. We put $\sigma_1^-:=\emptyset$ and for $i=2,3,\cdots,K$ the partition $\sigma_i^-$ is the part of $\sigma$ which lies to the left of $\Delta_i$, so $\sigma_i^-:=\bigcup\limits_{1\leq q < i}\sigma_q$. For $i=1,2,3,\cdots,K$ we denote by $\tau_i^+$ the part of $\tau$ which lies to the right of $\Delta_i$, so $\tau_i^+:=\bigcup\limits_{i< q \leq K}\Delta_q$ and $\tau_{K+1}^+:=\emptyset$. Finally for $i=1,2,\cdots,K$ we let $\tau(i):=\sigma_i^-\cup \Delta_i\cup \tau_i^+$ and $\tau(K+1):=\sigma$. 

Now we write the kernel $\Pi^{\tau(i)}$ as $\Pi_i$ and obtain in particular $\Pi_1=\Pi^\tau$, $\Pi_{K+1}=\Pi^\sigma$ and \[\Pi^\tau-\Pi^\sigma = \sum\limits_{i=1}\limits^{K} (\Pi_i-\Pi_{i+1}).\]
If we now interpret $\Pi^\emptyset$ as convolution identity operator we have for $i=1,2,\cdots,K$
\[\Pi_i-\Pi_{i+1}=\Pi^{\tau_{i}^+}\circ (\omega_{\Delta_i}-\Pi^{\sigma_i})\circ \Pi^{\sigma_i^-}.\]

If we now use \thref{Pi_mu} and \thref{pi_omega}, we obtain
\begin{align*}
|\Pi_i-\Pi_{i+1}|\leq & k_{|\Delta_i^+|}+\epsilon c(S) \frac{M(\Delta_i^+)^{\alpha-1}}{m(\Delta_i^+)^\alpha}|\Delta_i^+|k_{m(\Delta_i^+)}\\
&+\epsilon k_{m(\Delta_i^+)} |\Delta_i^+|^2\\
&(\epsilon \left(\frac{M(\Delta_i^+)+|\Delta_i^+|}{m(\Delta_i^+)}\right)^\alpha\sum\limits_{q=2}\limits^{\infty}\frac{q^\alpha}{q!}\left(\epsilon\frac{|\Delta_i^+|}{m(\Delta_i^+)}\right)^{q-2}\left(c(S)\frac{M(\Delta_i^+)^{\alpha-1}}{m(\Delta_i^+)^{\alpha-1}}\right)^q+\\
&+\epsilon c(S) \left(\frac{ M(\Delta_i^+)^{\alpha-1}(3M(\Delta_i^+)-m(\Delta_i^+))^\alpha}{m(\Delta_i^+)^{2\alpha+1}}+\frac{(M(\Delta_i^+)+|\Delta_i^+|)^{\alpha-1} }{m(\Delta_i^+)^{\alpha+1}}\right) )\\
&\circ \\
&\epsilon k_{m(\Delta_i)} |\Delta_i|^2\\
&(\epsilon \left(\frac{M(\Delta_i)+|\Delta_i|}{m(\Delta_i)}\right)^\alpha\sum\limits_{q=2}\limits^{\infty}\frac{q^\alpha}{q!}\left(\epsilon\frac{|\Delta_i|}{m(\Delta_i)}\right)^{q-2}\left(c(S)\frac{M(\Delta_i)^{\alpha-1}}{m(\Delta_i)^{\alpha-1}}\right)^q+\\
&+\epsilon c(S) \left(\frac{ M(\Delta_i)^{\alpha-1}(3M(\Delta_i)-m(\Delta_i))^\alpha}{m(\Delta_i)^{2\alpha+1}}+\frac{(M(\Delta_i)+|\Delta_i|)^{\alpha-1} }{m(\Delta_i)^{\alpha+1}}\right) ) \\
&\circ \\
&k_{|\Delta_i^-|}+\epsilon c(S) \frac{M(\Delta_i^-)^{\alpha-1}}{m(\Delta_i^-)^\alpha}|\Delta_i^-|k_{m(\Delta_i^-)}\\
&+\epsilon k_{m(\Delta_i^-)} |\Delta_i^-|^2\\
&(\epsilon \left(\frac{M(\Delta_i^-)+|\Delta_i^-|}{m(\Delta_i^-)}\right)^\alpha\sum\limits_{q=2}\limits^{\infty}\frac{q^\alpha}{q!}\left(\epsilon\frac{|\Delta_i^-|}{m(\Delta_i^-)}\right)^{q-2}\left(c(S)\frac{M(\Delta_i^-)^{\alpha-1}}{m(\Delta_i^-)^{\alpha-1}}\right)^q+\\
&+\epsilon c(S) \left(\frac{ M(\Delta_i^-)^{\alpha-1}(3M(\Delta_i^-)-m(\Delta_i^-))^\alpha}{m(\Delta_i^-)^{2\alpha+1}}+\frac{(M(\Delta_i^-)+|\Delta_i^-|)^{\alpha-1} }{m(\Delta_i^-)^{\alpha+1}}\right) ).\\
\end{align*}
The right hand side is bounded by 
\begin{align*}
|\Delta_i|^2 c(\Delta,S) k_{m(\Delta)}.
\end{align*}
Therefore
\begin{align*}
|\Pi^\tau-\Pi^\sigma|
&\leq  \sum\limits_{i=1}\limits^{K} |\Pi_i-\Pi_{i+1}|\\
&\leq \lambda(\tau)c(\Delta,S)k_{m(\Delta)}.
\end{align*}

\end{proof}
\end{lemma}

\paragraph{Existence of $\omega_\Delta$:}

We now define 
\begin{equation}
\label{limit}
\omega_\Delta:=\lim\limits_{n\rightarrow\infty}\Pi^{\tau_n(\Delta)}.
\end{equation} 
This limit exists as a uniform limit because 
\begin{equation*}
|\Pi^{\tau_n(\Delta)}-\Pi^{\tau_{n+1}(\Delta)}|\leq \frac{1}{2^n} c(\Delta , S)k_{m(\Delta)}
\end{equation*}
for all $n\in\mathbb{N}$ as in \thref{refinement} and 
\begin{equation*}
\lim\limits_{n\rightarrow\infty}\Pi^{\tau_n(\Delta)}=\Pi^{\tau_1(\Delta)}+(\Pi^{\tau_2(\Delta)}-\Pi^{\tau_1(\Delta)})+(\Pi^{\tau_3(\Delta)}-\Pi^{\tau_2(\Delta)})+\cdots.
\end{equation*}

\subsection{Properties of $\omega_\Delta$}
\begin{lemma}
\thlabel{omega_delta}
The kernel $\omega_\Delta$ has the following properties:
\begin{enumerate}
\item \label{o 1} $\omega_\Delta$ is continuous on $S\times S$
\item \label{o 2} $\int\limits_S \omega_\Delta (x,\zeta)d\omega^{z_0}(\zeta)=1$ for all $x\in S$
\item \label{o 3} for $0<a<b<c$ we have $\omega_{[a,c]}=\omega_{[b,c]}\circ \omega_{[a,b]}$
\item \label{o 4} $\omega_\Delta$ is positive for $|\Delta|>m(\Delta)$
\item \label{o 5} $|\omega_\Delta-\tilde{\omega}_\Delta|\leq c \epsilon^2 \left(\frac{|\Delta|}{m(\Delta)}\right)^2k_{m(\Delta)}$ for any segment $\Delta$ with $|\Delta|\leq m(\Delta)$
\end{enumerate}
\begin{proof}
Proof of \ref{o 1}:\\ As all $\tilde{\omega}_\Delta$ are continuous, $\omega_\Delta$ is continuous as a uniform limit of continuous functions. 

Proof of \ref{o 2}:\\ We observe that $\tilde{\Omega}_\Delta (1)=1$, the limit in (\ref{limit}) is uniform and $\omega^{z_0}$ is a probability measure. \\
Proof of \ref{o 3}:\\ We put $\Delta:=[a,c]$, $\Delta^-=[a,b]$, $\Delta^+=[b,c]$ and $\tau_n ' (\Delta):=\tau_n (\Delta^-)\cup \tau_n (\Delta^+)$. Note that $\tau_n ' (\Delta)$ is almost the same as $\tau_n(\Delta)$ with the exception of one segment that might be split into two by $b$. Therefore $\tau_n ' (\Delta)\succ \tau_n(\Delta)$, by \thref{refinement}
\begin{equation*}
|\Pi^{\tau_n ' (\Delta)}-\Pi^{\tau_n(\Delta)}|\leq c(S,\Delta)\frac{1}{2^n}k_{m(\Delta)}
\end{equation*}
and $\lim\limits_{n\rightarrow\infty} \Pi^{\tau_n ' (\Delta)}=\omega_{\Delta}$. \thref{Pi_mu} applied to $\tau_n (\Delta^+)$ and $\tau_n (\Delta^-)$ provides an integrable majorant so that the limit passage in 
\begin{equation*}
\lim\limits_{n\rightarrow\infty} \Pi^{\tau_n (\Delta^+)}\circ \Pi^{\tau_n (\Delta^-)}=\omega_{\Delta^+}\circ \omega_{\Delta^-}
\end{equation*}
is justified. Noting that $\Pi^{\tau_n (\Delta^+)}\circ \Pi^{\tau_n (\Delta^-)}=\Pi^{\tau_n' (\Delta)}$ We have shown $\omega_{[a,c]}=\omega_{[b,c]}\circ \omega_{[a,b]}$.\\
Proof of \ref{o 4}:\\ see section \ref{positivity}\\
Proof of \ref{o 5}:\\ By \thref{pi_omega} we know that there is a constant $c$  such that for any $\Delta$ with $|\Delta|\leq m(\Delta)$ and any regular partition $\mu$ we have $|\Pi^\mu-\tilde{\omega}_\Delta|\leq c \epsilon^2 \left(\frac{|\Delta|}{m(\Delta)}\right)^2k_{m(\Delta)}$. This is therefore valid for the dyadic partitions and  for the limit $\omega_\Delta$. 
\end{proof}
\end{lemma}

\subsection{Positivity of $\Pi^\mu$ and $\omega_\Delta$}\label{positivity}
First we prove that for any segment $\Delta$ with $m(\Delta)\leq|\Delta| \leq 3 m(\Delta)$ and any regular partition $\mu$ of $\Delta$ the kernel $\Pi^\mu$ is positive. 

As we already know we can rewrite $\Pi^\mu$:
\begin{equation*}
\Pi^\mu = k_{|\Delta|}
-\epsilon \sum\limits_{s=1}\limits^{K}  k_{M(\Delta)-M(j_s)} \circ b_{j_s}\circ k_{m(j_s)-m(\Delta)}
+ \sum\limits_{q=2}\limits^{K}\sum\limits_{l\in N_q}\pi_l.
\end{equation*}
For further calculations we will use the notation $v_s=k_{M(\Delta)-M(j_s)} \circ b_{j_s}\circ k_{m(j_s)-m(\Delta)}$ and $\rho_\mu=\sum\limits_{q=2}\limits^{K}\sum\limits_{l\in N_q}\pi_l$ and (\ref{pi}) is equivalent to

\begin{equation*}
\Pi^\mu =k_{|\Delta|}-\epsilon \sum\limits_{s=1}\limits^{K} v_s + \rho_\mu
\end{equation*}

As a first step we will show that $\sum\limits_{s=1}\limits^{K} v_s\leq c(S) k_{|\Delta|}$ and therefore analyse $v_s$ with the help of property \ref{bd 1} of $b_\Delta$:
\begin{align*}
|v_s|&= k_{M(\Delta)-M(j_s)} \circ |b_{<j_s}|\circ k_{m(j_s)-m(\Delta)}\\
&\leq k_{|\Delta|-|j_s|} \circ c(S) \frac{M(j_s)^{\alpha-1}}{m(j_s)^\alpha}|j_s|k_{m(j_s)}\\
&\leq k_{|\Delta|-|j_s|} \circ c(S) \frac{M(j_s)^{\alpha-1}}{m(\Delta)^\alpha}|j_s|k_{m(\Delta)}\\
&\leq c(S) \frac{M(\Delta)^{\alpha-1}}{m(\Delta)^\alpha}|j_s| k_{m(\Delta)+|\Delta|-|j_s|}\\
&\leq c(S)\frac{|j_s|}{m(\Delta)}k_{|\Delta|}.
\end{align*}
Summing up, we obtain
\begin{equation*}
\sum\limits_{s=1}\limits^{K} v_s\leq c(S) \frac{|\Delta|}{m(\Delta)} k_{|\Delta|}\leq c(S)k_{|\Delta|}.
\end{equation*}
The next step is the analysis of $\rho_\mu$. Here we use the regularity of the partition:
For $l\in N_q$ we use the following notation:
\begin{equation*}
a(l):=\sum\limits_{s\in l}m(j_s)+\sum\limits_{s\notin l}|j_s|
\end{equation*}
Now we obtain for $l\in N_q$
\begin{align*}
|\pi_l|&\leq k_{a(l)} \prod\limits_{s\in l}\left( \epsilon c(S) \frac{M(\Delta)^{\alpha-1}}{m(\Delta)^{\alpha-1}}\frac{|j_s|}{m(\Delta)} \right)\\
&\leq k_{a(l)}\frac{1}{K^q}\left(\epsilon c(S)\frac{|\Delta|}{m(\Delta)}\right)^q
\end{align*}
Using the abbreviation $\nu :=  c(S)\frac{|\Delta|}{m(\Delta)}$ and $R(\nu):=\sum\limits_{q=2}\limits^{\infty}\epsilon^{q-2}\frac{\nu ^q q^\alpha}{q!}$. This shows
\begin{equation*}
|\rho_\mu|\leq \epsilon^2\sum\limits_{q=2}\limits^{K} \frac{\epsilon^{q-2}\nu^q}{q!}\left(\frac{qM(\Delta)+|\Delta|}{m(\Delta)} \right)^\alpha k_m(\Delta)\leq \epsilon^2 c(S)R(\nu)k_m(\Delta)
\end{equation*}
and $R(\nu)$ decreases if $\epsilon$ decreases.

Proof that $\omega_\Delta>0$ for $|\Delta|>m(\Delta)$:
As $\Pi^\mu > 0$ for any segment $\tilde{\Delta}$ with $m(\tilde{\Delta)}\leq|\tilde{\Delta}| \leq 3 m(\tilde{\Delta)}$ we know that the corresponding limit $\omega_{\tilde{\Delta}}$ is positive. 
Therefore we partition our segment as follows:
$\Delta=[m(\Delta),2m(\Delta)]\cup [2m(\Delta),4m(\Delta)]\cup \cdots \cup [2^n m(\Delta),M(\Delta)] $ where $n$ is chosen such that $|[2^n m(\Delta),M(\Delta)]|\geq 2^n m(\Delta)$. 
Then we use property (3) of $\omega_\Delta$ to rewrite $\omega_\Delta=\omega_{[2^n m(\Delta),M(\Delta)]}\circ \cdots\circ \omega_{[2m(\Delta),4m(\Delta)]}\circ \omega_{[m(\Delta),2m(\Delta)]}$.
As all the involved kernels are positive $\omega_\Delta$ is also positive.

\subsection{Properties of $\Omega_\Delta$}
\begin{lemma}[$\Phi$-Property]
\thlabel{phi prop}
Let $\epsilon, y\in (0,1)$ and let $\psi$ be a function defined on $S$ coinciding with a positive harmonic function $v$ on $O_{-y}$ with $v\big|_{O}=K(\psi)$ (the harmonic extension of $\psi$ recovers $v$ on the near-half space). Then for any $\Delta \in \seg$ with $\Delta\subset(0,y]$ and $|\Delta|\leq m(\Delta)$ we have: 
\begin{equation*}
|\Omega_\Delta(\psi)-\psi|\leq c(S) \frac{\Delta}{y}\psi.
\end{equation*}
\begin{proof}
Let $J\in\seg$ with $J\subset \Delta$. Using the notation $\omega_\Delta=\tilde{\omega}_\Delta+r_\Delta$ we have
\begin{equation*}
|\Omega_J(\psi)-\psi|\leq |K_{|J|}(\psi)-\psi|+|\epsilon B_J(\psi)|+|R_J(\psi)|.
\end{equation*}
We now estimate the three terms separately, starting with $|K_{|J|}(\psi)-\psi|$:
\begin{equation*}
|K_{|J|}(\psi)-\psi|(x)=|v(x_{|J|})-v(x)|\leq |\nabla v(x_\eta)|\cdot|J|
\end{equation*}
where $\eta=\eta(x)\in (0,|J|)$. As $\dist(x_\eta, S_{-y})\geq c_1 y$ and $v$ is positive and harmonic on $O_{-y}$, by Harnack's inequality
\begin{equation*}
|\nabla v(x_\eta)|\leq c_2\frac{v(x_\eta)}{y}\leq c_3 \frac{v(x)}{y}=c_3\frac{\psi(x)}{y}. 
\end{equation*}
For the estimate of $|\epsilon B_J(\psi)|$ we first observe that for $\theta \in J$
\begin{equation*}
|C_\theta(\psi)(x)|\leq |\nabla v(x_\theta)|\leq c_4 \frac{v_\theta(x)}{y} =c_4 \frac{K_\theta(\psi)(x)}{y}.
\end{equation*}
By the definition of $b_\Delta$ (\ref{def b_delta}) we can estimate as follows: 
\begin{equation*}
|\epsilon B_J(\psi)|\leq \epsilon\int\limits_{J}K_\theta(|C_\theta(\psi)|)d\theta\leq \epsilon\frac{c_4}{y}\int\limits_{J} K_{2\theta} (\psi) d\theta=\epsilon\frac{c_4}{y}\int\limits_{J} v_{2\theta} d\theta\leq  \epsilon \frac{c_5}{y} v \cdot |J|.
\end{equation*}
Next, we will estimate $|R_J(\psi)|$ using \thref{omega_delta}(\ref{o 5}):
\begin{align*}
\begin{split}
|R_J(\psi)|(x)&\leq c_6 \epsilon^2 \frac{|J|^2}{m(J)^2} \int\limits_{S} k_{m(J)} (x,\zeta)\psi(\zeta)d\omega^{z_0}(\zeta)\\
&\leq c_6  \epsilon^2 \frac{|J||\Delta|}{m(\Delta)^2} v(x_{m(J)}) \leq c_7 \frac{|J|^2}{ m(\Delta)^2 y} \psi(x). 
\end{split}
\end{align*}

Collecting the estimates we get
\begin{equation}
\label{estimate J}
(1-\rho_J)\psi\leq \Omega_J(\psi)\leq (1+\rho_J)\psi
\end{equation}
with $\rho_J\leq c_8\frac{|J|}{y}(1+\frac{|J|}{m(\Delta)^2})$.

Let $K=K(\Delta,y)$ be a positive integer large enough such that $2c_8\frac{|\Delta|}{K y}=:\sigma_K<\frac{1}{2}$ and $\frac{|\Delta|}{K m(\Delta)^2}\leq 1$.
We will now decompose $\Delta$ into $K$ non-overlapping segments $J_1,\cdots,J_K$ of equal length such that
$\Delta=\bigcup \limits_{k=1}\limits^{K}J_k$, $|J_k|=\frac{|\Delta|}{K}$ and $m(J_k)<m(J_{k+1})$ for all $k$. 
Our choice of $K$ also implies that $\rho_{J_k}\leq \sigma_K$.

Now by \thref{omega_delta} (\ref{o 3}) we can decompose $\Omega_\Delta$ and use (\ref{estimate J}):
\begin{align*}
\begin{split}
\Omega_\Delta (\psi) &=\Omega_{J_K}\Omega_{J_{K-1}}\cdots\Omega _{J_1}(\psi)\\
&\leq (1+\sigma_K)^K \psi=(1+2 c_8\frac{|\Delta|}{Ky})^K\psi\\
&<e^{2 c_8\frac{|\Delta|}{y}}\psi<(1+c_9\frac{|\Delta|}{y}) \psi.
\end{split}
\end{align*}
Similarly, the estimate from below is given by
\begin{align*}
\begin{split}
\Omega_\Delta (\psi) &=\Omega_{J_K}\Omega_{J_{K-1}}\cdots\Omega _{J_1}(\psi)\\
&\geq (1-\sigma_K)^K \psi=(1-2 c_8\frac{|\Delta|}{Ky})^K\psi\\
&\geq e^{-c_{10}\frac{|\Delta|}{y}}\psi \geq(1-c_{10}\frac{|\Delta|}{y}) \psi.
\end{split}
\end{align*}
The estimates from above and below imply the assertion in the lemma. 
\end{proof}
\end{lemma}

\begin{lemma}
\thlabel{phi prop 2}
Let $\epsilon\in (0,\epsilon(S)$ and let $\psi\in\mathcal{C}(\bar{S})$. Then for any $\Delta \in \seg$ with $|\Delta|\leq m(\Delta)$ we have:
\begin{equation*}
\|\psi-\Omega_\Delta(\psi)\|_\infty\leq \|\psi-K_{|\Delta|}(\psi)\|_\infty+c(S)\frac{|\Delta|}{m(\Delta)}\|\psi\|_\infty.
\end{equation*} 

\begin{proof}
The proof is very similar to the one of \thref{phi prop}. We again start with
\begin{equation*}
|\Omega_\Delta(\psi)-\psi|\leq |K_{|\Delta|}(\psi)-\psi|+|\epsilon B_\Delta(\psi)|+|R_\Delta(\psi)|.
\end{equation*}
We do not estimate the first term. The second one is now done in the following way: 
\begin{equation*}
|\epsilon B_\Delta(\psi)|\leq \int\limits_{\Delta} |B_\theta|(|\psi|)\leq c(S)\int\limits_{\Delta} \frac{K_\theta(|\psi|)}{\theta} d\theta \leq c(S) \frac{|\Delta|}{m(\Delta)}\|\psi\|_\infty
\end{equation*}
For the third term we get (again by \thref{omega_delta}(\ref{o 5}))
\begin{equation*}
|R_\Delta(\psi)|(x)\leq c \epsilon^2 \frac{|\Delta|^2}{m(\Delta)^2} \int\limits_{S} k_{m(\Delta)} (x,\zeta)|\psi|(\zeta)d\omega^{z_0}(\zeta)\leq c\frac{\Delta}{m(\Delta)}\|\psi\|_\infty.
\end{equation*}
Similarly as in \thref{phi prop} we reach the conclusion by taking the supremum norm. 
\end{proof}
\end{lemma}

From now on we use the notation $\omega_y:=\omega_{[y,1]}$ and $\Omega_y:=\Omega_{[y,1]}$.
\begin{lemma}
\thlabel{omega_rho}
For $\rho\in (0,\frac{1}{2})$ and $\epsilon\in (0,\epsilon (S)$ we have the following estimates
\begin{align*}
\omega_\rho &\leq  c_+ \frac{1}{\rho^{c_+\epsilon}}k_{1-\rho}\\
\omega_\rho &\geq  c_-\rho ^{c_-\epsilon} k_{1-\rho}
\end{align*}
where $c_+(S)$ and $c_-(S)$ are positive constants. 
\begin{proof}
We start by estimating $\omega_{[y,2y]}$ for $y>0$. 
In the following estimates we use \thref{omega_delta} (\ref{o 5}), the definition of $\tilde{\omega}_\Delta$ and $b_\Delta$ and \thref{prop by} as well as \thref{prop ky}:
\begin{align}
\begin{split}
\label{absch fuer y 2y}
\omega_{[y,2y]}&\leq k_y+\epsilon |b_{[y,2y]}|+c\epsilon^2 k_y\\
&\leq (1+c\epsilon^2) k_y +\epsilon\int\limits_{y}\limits^{2y}\frac{k_\theta}{\theta}d\theta\\
&\leq (1+c\epsilon^2) k_y +\epsilon\int\limits_{y}\limits^{2y} k_y\frac{\theta^{\alpha-1}}{y^{\alpha}}d\theta\\
&\leq (1+c\epsilon^2+\epsilon (2^\alpha-1)) k_y \\
&\leq (1+c_2\epsilon)k_y
\end{split}
\end{align} 
Now let $K=K(\rho)$ be the natural number such that $2^K \rho\leq 1\leq 2^{K+1}\rho$. This implies $K\log (2)\leq \log(\frac{1}{\rho})\leq (K+1)\log (2)$.
We partition the segment $[\rho,1]$ as follows
\begin{equation*}
[\rho,1]=\left(\bigcup\limits_{j=0}\limits^{K-1} [2^j \rho, 2^{j+1} \rho]\right)\cup [2^K\rho, 1].
\end{equation*}
Using \thref{omega_delta}(\ref{o 3}) we obtain
\begin{equation*}
\omega_{[\rho, 1]}=\omega_{[2^K\rho, 1]}\circ \omega_{[2^{K-1}\rho,2^{K}\rho]}\circ \cdots \circ \omega_{[\rho, 2\rho]}
\end{equation*}
Now we use (\ref{absch fuer y 2y}) and obtain
\begin{equation*}
\omega_{[\rho, 1]}\leq (1+c\epsilon)^{K+1}k_{1-\rho}\leq 2\frac{1}{\rho^{c_1\epsilon}}k_{1-\rho}\leq c_+\frac{1}{\rho^{c_+\epsilon}}k_{1-\rho}
\end{equation*}

The second inequality can be proven similarly starting with $\omega_{[y,2y]}\leq k_y-\epsilon|b_{[y,2y]}|-c\epsilon^2 k_y$. 
\end{proof}
\end{lemma}

\begin{lemma}[Differential equation]
\thlabel{diff equ}
Let $\phi$ be a positive harmonic function with $\lim\limits_{z\rightarrow\infty}\phi(z)=0$. Then $f^x:y \mapsto \Omega_y(\phi_y)(x)$ is continuously differentiable and
\[\left(\frac{\partial}{\partial y}f^x\right)(y)=\epsilon \Omega_y(B_y(\phi_y))(x)\]
for $x\in S$ and $y\in (0,1]$.

\begin{proof}
We will start by proving that $f^x$ is Lipschitz on any segment $[y_0,1]$ for $0<y_0<1$. Then we will compute the left derivative  $(f^x)'_{-}$ which exist everywhere on $(0,1]$ and is continuous. By the Lipschitz property we obtain that $f^x(y)=f^x(1)-\int\limits_y\limits^{1}(f^x)'_{-}(\eta)d\eta$ for $y\in(0,1]$. Therefore the right and left derivative are equal and $f^x\in\mathcal{C}^1((0,1])$. 

We start by showing that $f^x$ is Lipschitz for $x\in S$ and $y_0\in(0,1]$. For $y\in (0,1]$, $h>0$, $y-h\geq y_0$ and $\Delta=[y-h,y]$ we have by definition of $\Omega_y$ and (\ref{o 3}):
\begin{equation*}
f^x(y)-f^x(y-h)=\Omega_y(I+II)(x)+III(x),
\end{equation*} 
where $I=\phi_y-\Omega_\Delta(\phi_y)$, $II=\phi_y-\phi_{y-h}$ and $III=(\Omega_{y-h}-\Omega_y)(II)$.

We set $\|\phi_y\|:=\sup\limits_{S}\phi_y$ and $K:=\sup \limits_{y_0\leq y\leq 1}\|\phi_y\|$ and we recall that by \thref{omega_delta} (\ref{o 2}) the norm of the operator $\Omega_y:\mathcal{C}([y,1])\rightarrow \mathcal{C}([y,1])$ does not exceed $1$.
Using the $\Phi$-Property in \thref{phi prop} we obtain 
\begin{equation*}
\|I\|\leq c_1(S)K\frac{h}{y_0}.
\end{equation*}
For any $x\in S$ by the mean value theorem there is a $\theta=\theta(x)\in (y-h,y)$ such that
\begin{equation*}
|II|(x)\leq |\nabla\phi(x_\theta)|h\leq c_2(S)K\frac{h}{y_0},
\end{equation*}
where the second inequality holds because of Harnack's inequality. Hence
\begin{equation*}
\|\Omega_y(I+II)\|\leq \|I\|+\|II\|\leq c_3(S)K\frac{h}{y_0}.
\end{equation*}
Finally, 
\begin{equation*}
\|III\|\leq 2\|II\|\leq 2c_2(S)K\frac{h}{y_0}
\end{equation*}
and therefore $f^x$ is Lipschitz on $[y_0,1]$ for $x\in S$ because
\begin{equation*}
|f^x(y)-f^x(y-h)|\leq c_4(S)K\frac{h}{y_0}.
\end{equation*}

Now we will compute the left derivative of $f^x$. For $x\in S$, $y\in (0,1)$ and $h\in (0,\frac{y}{2})$ we have
\begin{equation*}
\frac{f^x(y-h)-f^x(y)}{-h}=\Omega_y\left(\frac{I+II}{h}\right)(x)+\frac{III}{h}(x).
\end{equation*}
As a first step we will prove that $\lim\limits_{h\downarrow 0}\frac{III}{h}(x)=0$ on $S$. In order to do this we rewrite $\frac{III}{h}=(\Omega_y-\Omega_{y}\Omega_\Delta)(\frac{\partial \phi}{\partial \vec{e}_d}+\frac{\phi_y-\phi_{y-h}}{h}-\frac{\partial \phi}{\partial \vec{e}_d})$. 
As $\frac{\partial \phi}{\partial \vec{e}_d}\in\mathcal{C}(\bar{S})$, 
$|\frac{\partial \phi}{\partial \vec{e}_d}(S_y)|\leq c\frac{\phi|_{S_y}}{y}$, 
$\lim\limits_{\infty}\frac{\partial \phi}{\partial \vec{e}_d}\big|_{S_y}=0$ we can apply \thref{phi prop 2} to obtain a uniform limit $0$ as $h\downarrow 0$ of $\frac{\partial \phi}{\partial \vec{e}_d}-\Omega_{\Delta}(\frac{\partial \phi}{\partial \vec{e}_d})$. Therefore $(\Omega_y-\Omega_{y}\Omega_\Delta)(\frac{\partial \phi}{\partial \vec{e}_d})$ converges to $0$ as $h\downarrow 0$. 
Also, since $\phi$ is bounded on $O_{\frac{y}{2}}$ and its second derivatives are bounded on any ball of radius $\frac{y}{4}$ and center on $S+\vec{e}_d$ by some constant only depending on $y$, we know $|\frac{\phi_y-\phi_{y-h}}{h}-\frac{\partial \phi}{\partial \vec{e}_d}|\leq c(y)h$ on $S$. Therefore $|(\Omega_y-\Omega_{y}\Omega_\Delta)(\frac{\phi_y-\phi_{y-h}}{h}-\frac{\partial \phi}{\partial \vec{e}_d})|\leq 2 c(y) h$ and we proved  $\lim\limits_{h\downarrow 0}\frac{III}{h}(x)=0$ on $S$.

As a second step, we will prove that $\lim\limits_{h\downarrow 0} \frac{I+II}{h}=\epsilon B_y(\phi_y)$ on $S$. Using the notation $\omega_\Delta=\tilde{\omega}_\Delta +r_\Delta$, we have 
\begin{align*}
\begin{split}
\frac{I+II}{h}&=\frac{\phi_y-\Omega_\Delta(\phi_y) +\phi_y-\phi_{y-h}}{h}\\
&=\frac{\phi_y-K_h(\phi_y)+\epsilon B_\Delta(\phi_y)-R_\Delta(\phi_y) +\phi_y-\phi_{y-h}}{h}.
\end{split}
\end{align*}
As $K_h(\phi_y)=\phi_{y+h}$ we have $\lim\limits_{h\downarrow 0} \frac{\phi_y-K_h(\phi_y)+\phi_y-\phi_{y-h}}{h}=0$ and using the definition of $b_\Delta$ (\ref{def b_delta})and the continuity of $b_y$ \thref{prop by} (\ref{by 3}) we obtain
\begin{equation*}
\lim\limits_{h\downarrow 0} \frac{I+II}{h}=\lim\limits_{h\downarrow 0}\frac{\epsilon}{h}\int \limits_{y-h}\limits^{y} B_\theta(\phi_y)d\theta + \frac{1}{h}R_\Delta(\phi_y)=B_y(\phi_y)+ \lim\limits_{h\downarrow 0}\frac{1}{h}R_\Delta(\phi_y).
\end{equation*}
The last limit is zero by \thref{omega_delta} (\ref{o 5}) and the calculation of the left derivative and therefore the proof is finished.

\end{proof}
\end{lemma}

\begin{lemma}
Let $\phi$ be a positive harmonic function on $O$ with a finite limit $\lim\limits_{z\rightarrow\infty}\phi(z)$. Then for $0<\eta<y\leq 1$
\[\Omega_\eta(\phi_y)\leq c(S)\Omega_y(\phi_y)\]
\begin{proof}
This is a direct consequence of the $\Phi$-property \thref{phi prop}. We put $\Delta=[\eta,y]$ and $\psi=\phi_y$ so that 
\begin{equation*}
\Omega_\eta(\phi_y)=\Omega_y(\Omega_\Delta(\phi_y))\leq (1+c(S))\Omega_y(\phi_y).
\end{equation*}

\end{proof}
\end{lemma}

\subsection{Measures}
For any probability measure $\kappa$ on $S$ and a fixed $\epsilon\in (0,\epsilon (S))$ we obtain a transformed measure as follows:
\[\gamma_y(x):=\int\limits_{S}\omega_y(\zeta,x)d\kappa(\zeta)=\Omega_y^*(\kappa)(x);\]
the measures we are interested in are the ones with density $\gamma_y$ with respect to the harmonic measure $\omega^{z_0}$.
The limit of those will be our desired measure $\nu_\epsilon$:
\[d\nu_\epsilon=\lim\limits_{y\downarrow 0 } \gamma_y d\omega^{z_0}.\]

\paragraph*{Existence of $\nu_\epsilon$}

We will now prove the weak convergence of the measures with density $\gamma_y$.
\begin{lemma}
The measures with density $\gamma_y$ with respect to the harmonic measure $\omega^{z_0}$ converge weakly as $y \downarrow 0$ to some measure $\nu_\epsilon$ on $\bar{S}$ supported on $S$ and $\nu_\epsilon(S)=1$.  
\begin{proof}
First we choose a monotone decreasing sequence $(y_k)_{k\in \mathbb{N}}$ in $(0,1)$ which converges to zero and such that the measures with density $\gamma_{y_k}$ w.r.t the harmonic measure converge weakly to some measure $\nu_\epsilon$ on $\bar{S}$, so 
\begin{equation*}
\lim\limits_{k\rightarrow \infty} \int\limits_{S} \alpha \gamma_{y_k} d\omega^{z_0} =\int\limits_{\bar{S}} \alpha d\nu_\epsilon
\end{equation*}
for all continuous functions $\alpha$. We will now verify that $\nu_\epsilon(\{\infty\})=0$ so that we can write $S$ instead of $\bar{S}$ in the integral on the right-hand side and $\nu_\epsilon$ is a probability measure on $S$. 
We consider a ball $\mathbb{B}_L$ with radius $L$ large enough so that $S$ is flat outside of $\mathbb{B}_L$ (so $S\setminus\mathbb{B}_L\subset \mathbb{R}^{d-1}$). This is possible because of the definition of the geometry of our near-half space $O$. 
We now consider the harmonic measures $\omega^z(\mathbb{B}_L\cap S,O)$ and $\omega^z(S\setminus\mathbb{B}_L,O)$ as harmonic functions of $z$ in $O$. We note that the sum of those two functions is always $1$. As $\omega^z(\mathbb{B}_L\cap S,O)$ vanishes on $S\setminus\bar{\mathbb{B}}_L$ it admits a harmonic extension to the domain $O\cup (\mathbb{R}^d\setminus \bar{\mathbb{B}}_L)$. The extension is bounded and vanishes at infinity. We now choose $\rho>0$ large enough so that $\{x_d=\rho\}\subset O$ and $L'>L$ so large that $(\omega^x(S\setminus\mathbb{B}_L,O))_\rho>\frac{1}{2}$ for all $x\in S\setminus \bar{\mathbb{B}}_{L'}=\mathbb{R}^{d-1}\setminus \bar{\mathbb{B}}_{L'}$. Putting $\omega^\infty(S\setminus\mathbb{B}_L,O)=1$ we may assume that $(\omega^\cdot(S\setminus\mathbb{B}_L,O))_\rho \big|_{\bar{S}}\in \mathcal{C}(\bar{S})$.
We now have 
\begin{equation*}
\nu_\epsilon(\{\infty\})\leq \nu_\epsilon(\bar{S}\setminus \mathbb{B}_{L'})\leq 2\int\limits_{\bar{S}}(\omega^\cdot(S\setminus\mathbb{B}_L,O))_\rho d\nu_\epsilon.
\end{equation*}
By definition of $\nu_\epsilon$ and $\gamma_y$ this is 
\begin{align*}
 2\int\limits_{\bar{S}}(\omega^\cdot(S\setminus\mathbb{B}_L,O))_\rho d\nu_\epsilon
 &=\lim\limits_{k\rightarrow\infty} \int\limits_{S} (\omega^\cdot(S\setminus\mathbb{B}_L,O))_\rho\gamma_{y_k} d\omega^{z_0}\\
 &=\lim\limits_{k\rightarrow\infty} \int\limits_{S}\Omega_{y_k}( (\omega^\cdot(S\setminus\mathbb{B}_L,O))_\rho) d\kappa
\end{align*}
Now we use \thref{phi prop} to obtain the following estimate for $\nu_\epsilon(\{\infty\})$
\begin{equation*}
\lim\limits_{k\rightarrow\infty} \int\limits_{S}\Omega_{y_k}( (\omega^\cdot(S\setminus\mathbb{B}_L,O))_\rho) d\kappa\leq \frac{c(S)}{\rho}\int\limits_{S} \omega^{x+\rho \vec{e}_d} (S\setminus\mathbb{B}_L,O)d\kappa(x).
\end{equation*}
For $L\rightarrow +\infty$ the harmonic measure of $S\setminus\mathbb{B}_L$ tends to $0$ on $S_\rho$ and is bounded by $0$ from below and $1$ from above. Considering that $\kappa$ is a probability measure on $S$ the estimate for $\nu_\epsilon(\{\infty\})$ tends to $0$ as $L\rightarrow +\infty$ and $\nu_\epsilon(\{\infty\})$ is zero.

It remains to show that 
\begin{equation}
\label{alpha}
\lim\limits_{y\downarrow 0}\int \limits_{S}\alpha\gamma_yd\omega^{z_0} =\int\limits_{S} \alpha d\nu_\epsilon
\end{equation}
for all $\alpha\in \mathcal{C}(\bar{S})$. First we assume that $\alpha$ coincides with $K_\sigma(\psi)$ for some positive $\psi \in \mathcal{C}(\bar{S})$ and $\sigma>0$. We can exploit \thref{phi prop}, \thref{omega_delta}(\ref{o 3},\ref{o 2}) and the fact that $\kappa$ is a probability measure to show that for all $0<y<y_k<1/2$
\begin{align*}
\begin{split}
&\left|\int\limits_{S}\alpha \gamma_y d\omega^{z_0} - \int\limits_{S}\alpha \gamma_{y_k} d\omega^{z_0} \right|=\\
&= \left|\int\limits_{S}\Omega_{y_k}(\Omega_{[y,y_k]}(\alpha)-\alpha)d\kappa \right|\leq \\
&\leq \|\Omega_{[y,y_k]}(\alpha)-\alpha\|_{\infty,S}\leq c(S)\frac{y_k}{\sigma}\|\alpha\|_{\infty,S}.
\end{split}
\end{align*}
As $y_k\rightarrow 0$ as $k\rightarrow \infty$, the convergence is proven in this case. It remains to note that for any $\alpha\in \mathcal{C}(\bar{S})$ we have $\|K_\sigma(\alpha)-\alpha\|_{\infty,S}\rightarrow 0$ as $\sigma\downarrow 0$.

\end{proof}
\end{lemma}

\paragraph*{Properties of $\nu_\epsilon$:}
\begin{lemma} The measure $\nu_\epsilon$ has the following properties:
\thlabel{measure}
\begin{enumerate}

\item \label{nu 1} For any $\epsilon\in (0,\epsilon(S))$ we have \[\int\limits_{S} V(x) d\nu_\epsilon(x)\leq\frac{c}{\epsilon}\int\limits_{S} u_1 d\kappa\]

\item \label{nu 3} For any ball $\mathbb{B}$ with center on $S$ there is an $\epsilon(\mathbb{B})$ such that \[\nu_\epsilon(\mathbb{B})>c\] for any $0<\epsilon<\epsilon(\mathbb{B})$. The constant $c$ may depend on the radius of the ball, $S$ and the probability measure $\kappa$. 

\end{enumerate}

\begin{proof}
Proof of \ref{nu 1}
For $y\in (0,1]$ we put \[g_y:=B_y(u_y).\] For $\delta\in (0,1)$ we want to prove that \[J_\delta:=\int\limits_{S}\int\limits_{\delta}\limits^{1}g_y dy d\nu_\epsilon\] is uniformly bounded. We obtain:
\begin{align*}
J_\delta&=\lim\limits_{\eta\rightarrow 0} \int\limits_{S}\int\limits_{\delta}\limits^{1}g_y dy \gamma_\eta d\omega^{z_0}=\lim\limits_{\eta\rightarrow 0} \int\limits_{S} \Omega_\eta\left( \int\limits_{\delta}\limits^{1}g_y dy \right) d\kappa=\\
&=\lim\limits_{\eta\rightarrow 0} \int\limits_{S}  \int\limits_{\delta}\limits^{1}\Omega_\eta\left(g_y\right) dy  d\kappa\leq 
c \int\limits_{S}  \int\limits_{\delta}\limits^{1}\Omega_y\left(g_y\right) dy  d\kappa=\\
&=\frac{c}{\epsilon}\int\limits_{S}\int\limits_{\delta}\limits^{1} \frac{\partial}{\partial y}\Omega_y(u_y)dy d\kappa=\frac{c}{\epsilon}\int\limits_{S}\Omega_1(u_1)-\Omega_\delta(u_\delta) d\kappa\leq\\
&\leq \frac{c}{\epsilon}\int\limits_{S}\Omega_1(u_1)d\kappa.
\end{align*}

Proof of \ref{nu 3}: 
Let $\zeta\in S$ and $r<1/2$ be the center and radius of the ball $\mathbb{B}$. Let $\psi$ be a function with $0\leq \psi\leq 1$, $\psi\equiv 1$ on $\frac{1}{2}\mathbb{B}$, $\psi \equiv 0$ outside of $\mathbb{B}$ and $|\nabla \psi|\leq \frac{2}{r}$. Let $\phi$ be a function on $S$ coinciding with $\psi$. As usual, we will denote by $\phi$ also the harmonic extension of $\phi$ to $O$.

For $y\in (0,r)$ we will now consider
\begin{align*}
\Omega_y(\phi_y)&=\Omega_r(\phi_r)-\int\limits_{y}\limits^{r}\frac{d}{d\theta}(\Omega_\theta(\phi_\theta))d\theta\\
&=\Omega_r(\phi_r)-\epsilon\int\limits_{y}\limits^{r}\Omega_\theta(B_\theta(\phi_\theta))d\theta,
\end{align*}
where we used the differential equation in \thref{diff equ}.
We will now estimate $\Omega_r(\phi_r)$ from below and $\epsilon\int\limits_{y}\limits^{r}\Omega_\theta(B_\theta(\phi_\theta))d\theta$ from above. 
By \thref{omega_rho}, we have
\begin{equation*}
\Omega_r(\phi_r)\geq c(S)r^\epsilon K_{1-r}(\phi_r)=c(S)r^\epsilon \phi_1
\end{equation*}
and therefore
\begin{equation*}
\int\limits_{S}\Omega_r(\phi_r)d\kappa\geq c_2 r^\epsilon.
\end{equation*}
On the other hand, 
\begin{equation*}
|\epsilon\int\limits_{y}\limits^{r}\Omega_\theta(B_\theta(\phi_\theta))d\theta|\leq \epsilon\int\limits_{y}\limits^{r}\Omega_\theta(|B_\theta(\phi_\theta)|)d\theta\leq \epsilon r \sup\limits_{S, \theta\in (0,r)} |B_\theta(\phi_\theta)|.
\end{equation*}
Here we used the positivity of $\Omega_\theta$ and the fact that $\Omega_\theta(1)=1$.
Estimating $|B_\theta(\phi_\theta)|$ using the definition of the kernels $k_y, c_y$ and $b_y$ and their properties, we obtain
\begin{equation*}
|B_\theta(\phi_\theta)|\leq c_3(S) \sup \limits_{O}|\nabla \phi|\leq \frac{1}{r}c_4(S,\mathbb{B}).
\end{equation*}

Collecting the estimates we obtain 
\begin{equation*}
\int \limits_{S}\Omega_y(\phi_y)d\kappa\geq c_2 r^\epsilon+c_4\epsilon
\end{equation*}
which is larger than a constant for $\epsilon\in (0,\epsilon(\mathbb{B},\kappa)$.

By \thref{phi prop} we know that 
\begin{align*}
\int \limits_{S}\Omega_y(\phi_y)d\kappa&= \int\limits_{S}\phi_y \Omega^*_y(\kappa)d\omega^{z_0}\leq \\
&\leq c(S)\int\limits_{S}\Omega_{[\delta,y]}\phi_y \Omega^*_y(\kappa)d\omega^{z_0}=\\
&=c(S)\int\limits_{S}\phi_y \Omega^*_\delta(\kappa)d\omega^{z_0}
\end{align*}

As we know that the first term is larger than a constant, this applies also to the last and as the limit of the measures $\Omega^*_\delta(\kappa)$ is the density of $\nu_\epsilon$ we obtain 
\begin{equation*}
\int\limits_{S}\phi_y d\nu_\epsilon>c
\end{equation*}
for all $y\in (0,r)$ and therefore the proof is finished. 

\end{proof}
\end{lemma}

\subsection{Proof of the Main Lemma}
\label{Proof of Lemma}
We now use \thref{measure} to prove \thref{main lemma}.
As the probability measure $\kappa$ in \thref{measure} was arbitrary, we can choose the harmonic measure $\omega^{z_1-\vec{e}_d}$ and obtain 
\[\int\limits_{S} V(x) d\nu_\epsilon(x)\leq c u(z_1).\]
As $\nu_\epsilon(\mathbb{B})>c(r(\mathbb{B})$ there is a point $x\in S$ such that \[V(x)\leq c u(z_1)\] where the constant $c$ my depend on the Lipschitz constant of the function defining $O$ and the radius $r(\mathbb{B})$ of the ball.

\subsection{Proof of Theorem 2}
Given a point $p_0$ on the boundary we choose a radius $r_0$ and a height $h_0$ such that (up to translation and rotation of the space) the cylinder with center $p_0=0$ radius $r_0$ and height $h_0$ \[\mathcal{C}(p_0):=\{(x,y)\in \mathbb{R}^{d-1}\times \mathbb{R}: |x|<r_0, |y|<h_0\}\] and the Lipschitz function $\phi$ ($\phi(0)=0$, $|\phi(x)|<\frac{h_0}{2}$) have the following properties:
\begin{itemize}
\item $\mathcal{C}^+:=\mathcal{C}(p_0)\cap O= \{(x,y)\in \mathbb{R}^{d-1}\times \mathbb{R}: |x|<r_0, \phi(x)<y<h_0\}$
\item $\mathcal{C}^-:=\{(x,y)\in \mathbb{R}^{d-1}\times \mathbb{R}: |x|<r_0,-h_0<y< \phi(x)\}\subset\mathbb{R}^d\setminus \bar{O}$
\item $S\cap \mathcal{C}$ is the graph of $\phi|_{\mathbb{B}(0,r_0)}$.
\end{itemize}
The argument used for the near-half space is now also applicable to the upper part of the cylinder $\mathcal{C}^+$ as it is a Lipschitz domain and because of its special geometry. In particular, if we take the cylinder with half the height of the original one, we can shift it in direction of $\vec{e}_d$ and remain inside $O$, where the function $u$ is positive and harmonic.

\subsection*{Acknowledgements}
This paper is part of the second named author’s PhD thesis written at the Department of Analysis, Johannes Kepler University Linz. The  research has  been  supported  by  the  Austrian  Science  foundation  (FWF) Pr.Nr P28352-N32.

\bibliographystyle{abbrv}

\end{document}